\documentclass[12pt]{amsart}
\begin{document}
\renewcommand{\Box}{\qed}

\def\Ave{\operatorname{Ave}}

\title{LOW $M^{\ast }$-ESTIMATES ON COORDINATE SUBSPACES}

\author{A. A. Giannopoulos}
\address{A. A. Giannopoulos, Department of Mathematics, University of Crete,
Iraklion, Crete, Greece}
\email{deligia@talos.cc.uch.gr}
\author{V. D. Milman}
\address{V. D. Milman, 
Department of Mathematics, Tel Aviv University, Tel Aviv, Israel}
\email{vitali@math.tau.ac.il}
\curraddr{(both authors)\vbox to 24pt{} Mathematical Sciences Research Institute,
1000 Centennial Road, Berkeley, CA 94720}
\thanks{This work was initiated
at the Institute for Advanced Study and completed
at the Mathematical Sciences Research Institute.
Research at MSRI is supported in part by NSF grant DMS-9022140.}

\begin{abstract}
Let $K$ be a symmetric convex body in ${\mathbf R}^n$. It is well-known
that for every $\theta\in (0,1)$
there exists a subspace $F$ of ${\mathbf R}^n$ with $\dim F=
[(1-\theta )n]$ such that 
$${\mathcal P}_F(K)\supseteq \frac{c\sqrt{\theta }}
{M_K}D_n\cap F,\leqno (\ast )$$
\noindent where ${\mathcal P}_F$ denotes the orthogonal
projection onto $F$. Consider a fixed coordinate system
in ${\mathbf R}^n$. We study the question whether an analogue
of ($\ast $) can be obtained when one is restricted to choose
$F$ among the coordinate subspaces ${\mathbf R}^{\sigma },\;
\sigma\subseteq\{1,\ldots,n\}$, with $|\sigma |=[(1-\theta )n]$.
We prove several ``coordinate versions" of ($\ast $) in terms
of the cotype-2 constant, of the volume ratio and other
parameters of $K$. The basic source of our estimates is an
exact coordinate analogue of ($\ast $) in the ellipsoidal
case. Applications to the computation of the number of
lattice points inside a convex body are considered throughout
the paper.
\end{abstract}

\maketitle

\section{Introduction}

{\it Notation}. Our setting is ${\mathbf R}^n$ equipped with an
inner product $\langle .,.\rangle $ and the associated Euclidean norm
defined by $|x|=\langle x,x\rangle^{1/2},x\in {\mathbf R}^n$. We denote
the Euclidean unit ball and the unit sphere by $D_n$ and $S^{n-1}$
respectively, and we write $\sigma $ for the rotationally invariant
probability measure on $S^{n-1}$. 

Let $K$ be a symmetric convex body in ${\mathbf R}^n$. Then, $K$ induces in a
natural way a norm $\|.\|_K$ on ${\mathbf R}^n$. In what follows we shall denote by $X_K$ the
normed space $({\mathbf R}^n,\|.\|_K)$.
As usual, $K^o=\{y\in {\mathbf R}^n:\langle y,x\rangle\leq 1 \;{\rm
for \;every} \;x\in K\}$ is the polar body of $K$, and $X_{K^o}=
({\mathbf R}^n,\|.\|_{K^o})$ is the dual space of $X_K$.

Finally, we consider the integral parameters
$$M=M_K=\left(\int_{S^{n-1}}\|x\|_K^2\;\sigma (dx)\right)^{1/2}\;\;\;\;,\;\;\;\;
M^{\ast }=M_{K^o}=\left(\int_{S^{n-1}}\|x\|_{K^o}^2\;\sigma (dx)\right)^{1/2},$$
\noindent which are up to a constant the mean widths of $K^o$ and $K$ 
respectively.

\bigskip

{\it Results}. The following inequality of the second named author
plays an important role in developing
a proportional
theory of high-dimensional convex bodies:

\noindent {\bf Theorem A} ({\it Low $M^{\ast }$-estimate}). {\it There
exists a function $f:(0,1)\rightarrow {\mathbf R}^+$ such that for every 
symmetric convex body $K$ in ${\mathbf R}^n$ and for every $\theta\in (0,1)$
one can find a subspace $F$ of ${\mathbf R}^n$ with $\dim F= [(1-\theta )n]$
satisfying}
$$\|x\|_K\geq\frac{f({\theta})}{M_{K^o}}\;|x|\;\;\;,\;\;\;x\in F.\leqno (1.1)$$

Theorem A was originally proved in [M1] and a second proof using the
isoperimetric inequality on $S^{n-1}$ was given in [M2] where it
was shown that (1.1) holds with $f(\theta )\geq c\theta $ for
some absolute constant $c>0$ (and with an estimate $f(\theta )\geq\theta +o(1-\theta )$ as
$\theta\rightarrow 1^-$). This was later improved to $f(\theta )\geq c\sqrt{\theta }$ in [PT], see also [M3] for a different proof with this best possible
$\sqrt{\theta }$-dependence. Finally, it was proved in [Go] that one
can have
$$f(\theta )\geq\sqrt{\theta }(1+O(\frac{1}{\theta n})).\leqno (1.2)$$

Moreover, if we fix some $\theta \in (0,1)$ and consider the
Grassmannian manifold $G_{n,k}$  of all $k$-dimensional subspaces
of ${\mathbf R}^n$, where $k=k(\theta ,n)=[(1-\theta )n]$, equipped with the Haar probability
measure $\nu_{n,k} $, then (1.1) holds true with $f(\theta )\geq c\sqrt{\theta }$
for all subspaces $F$ in a subset ${\mathcal A}_{n,k}$ of $G_{n,k}$ which is of
almost full measure $\nu_{n,k} ({\mathcal A}_{n,k})$ as $n\rightarrow\infty $.

Interchanging the roles of $K$ and $K^o$, we may equivalently read
Theorem A in the following geometric form:
$${\mathcal P}_F(K)\supseteq\frac{c\sqrt{\theta }}{M_K}\;D_n\cap F,\leqno (1.3)$$
\noindent where ${\mathcal P}_F$ denotes the orthogonal projection onto $F$. In
this paper we will follow the tradition and continue calling an inclusion of the type (1.3) a ``low $M^{\ast }$-estimate" (for $K^o$).

Among other applications of (1.3), let us mention the quotient of
subspace theorem and the reverse Santal\'{o} inequality [M1], [BM].

\bigskip

Let $\{e_1,\ldots,e_n\}$ be an arbitrary but
fixed orthonormal basis of ${\mathbf R}^n$
with respect to $\langle .,.\rangle $. For a subset $\sigma\subseteq
\{1,\ldots,n\}$ we naturally define the coordinate subspace ${\mathbf R}^{\sigma }=
\{x\in {\mathbf R}^n:\langle x,e_j\rangle =0\;{\rm if}\;j\notin\sigma \}$.
We write $D_{\sigma }$ for $D_n\cap {\mathbf R}^{\sigma }$ and
$Q_{\sigma }$ for the unit cube $Q_n\cap {\mathbf R}^{\sigma }=[-1,1]^
{\sigma }$ in ${\mathbf R}^{\sigma }$.

Our purpose is to discuss ``low $M^{\ast }$-estimates'' in the form
(1.3) when one is
restricted to choose $F$ among the coordinate subspaces of ${\mathbf R}^n$
of a certain dimension $m$ proportional to $n$.

\medskip 

In Section 2 we study the case of an ellipsoid $E$ in ${\mathbf R}^n$. It turns out that for {\it any} orthonormal basis of ${\mathbf R}^n$ one has results
analogous to (1.3) with almost the same $\sqrt{\theta }$-dependence
on the parameter $\theta $:

\medskip

\noindent {\bf Theorem B} ({\it Coordinate low $M^{\ast }$-estimate for
ellipsoids}).
{\it Let $E$ be an ellipsoid in
${\mathbf R}^n$ and $\theta\in (0,1)$. Then, there exists $\sigma\subseteq\{1,\ldots,n\},\;|\sigma |\geq (1-\theta )n$, with
$${\mathcal P}_{\sigma }(E)\supseteq\frac{c\sqrt{\theta }}{\log^{1/2}(\frac{2}{\theta })
M_E}\;D_{\sigma },$$
\noindent where ${\mathcal P}_{\sigma }$ denotes the orthogonal projection onto
${\mathbf R}^{\sigma }$, and $c>0$ is an absolute constant}.

\medskip

It is perhaps surprising that this type of geometric result about
ellipsoids is new and non-trivial. Note that our investigation
of these questions was started from a simpler fact of the same
nature about a special class of ellipsoids, which was discovered
in [Gi].

It can be checked that Theorem B is optimal apart from the logarithmic
term (see Remark 2.5). A result of the same type can be
proved for an ellipsoid $E$
of smaller but sufficiently large dimension living in an arbitrary
subspace $F$ of ${\mathbf R}^n$ (Theorem 2.3). We also consider the corresponding
problem for sections (instead of projections) of $E$ with coordinate subspaces (Theorem 2.4).

\medskip

Simple examples show that one cannot achieve the same strong estimate
in full generality: for an arbitrary symmetric convex body $K$ and
an arbitrary orthonormal basis in ${\mathbf R}^n$. Consider e.g the case
of the unit cube $Q_n$ and the standard basis of ${\mathbf R}^n$: observe
that $M_{Q_n}\simeq \sqrt{\log n/n}$, while the radius of the largest
Euclidean ball contained in any coordinate projection of $Q_n$ is 1.
In Section 3 we give a general low $M^{\ast }$-estimate in terms of
the cotype-2 constant $C_K$ of $X_K$:

\medskip

\noindent {\bf Theorem C} ({\it $M^{\ast }$-estimate in terms of $C_K$}). {\it For an arbitrary symmetric convex body $K$ in ${\mathbf R}^n$
and for any $\theta\in (0,1)$, one can find $\sigma\subseteq\{1,\ldots ,n\}$,
$|\sigma |\geq (1-\theta )n$, satisfying
$${\mathcal P}_{\sigma }(K)\supseteq\frac{c_1\theta }{\log^2 (\frac{2}{\theta })h(C_K)M_K}\;D_{\sigma },$$
\noindent where $h(y)=y\log 2y,\;y\geq 1$, and $c_1>0$ is an absolute
constant}.

\medskip

Let us note that one can give a simpler argument, based on the
isomorphic Sauer-Shelah lemma of S. J. Szarek and M. Talagrand
and a factorization theorem of B. Maurey, which results in a
weaker estimate of the same type (we sketch it in Remark 3.6).
We also obtain results of the same nature in which $M_K$ is
replaced by various other ``volumic" parameters of $K$ or $K^o$
(see Remark 3.7).

\medskip

In Section 4 we give a general low $M^{\ast }$-estimate in terms
of the volume ratio $vr(K)$ of $K$:

\medskip

\noindent {\bf Theorem D} ({\it $M^{\ast }$-estimate in terms of
$vr(K)$}). {\it Let $K$ be a symmetric convex body in ${\mathbf R}^n$.
For every $\theta\in (0,1)$, there exists $\sigma\subseteq\{1,\ldots,n\},
\;|\sigma |\geq (1-\theta )n$, such that
$${\mathcal P}_{\sigma }(K)\supseteq\frac{1}{[c_2vr(K)]^{\frac{c_3\log(\frac{2}{\theta })}{\theta }}M_K}\;D_{\sigma },$$
\noindent where $c_2,c_3>0$ are absolute constants}.

\medskip

In Sections 5 and 6 we give some further applications of the low
$M^{\ast }$-estimate for ellipsoids. We demonstrate an exact
dependence between coordinate sections of an ellipsoid and its
polar in the spirit of [M5]. We also apply Theorems 2.2 and 2.4
to questions related to the number of integer or ``almost integer"
points inside an ellipsoid.

\medskip 

Recall that the cotype-2 constant $C_K$ of $X_K$ is the smallest constant $\lambda >0$ for which
$$\left({\rm Ave }_{\varepsilon_j=\pm 1}\|\sum_{j=1}^m\varepsilon_jx_j\|_{K}^2
\right)^{1/2}\geq\frac{1}{\lambda }\left(\sum_{j=1}^m\|x_j\|_K^2\right)^{1/2}$$
\noindent holds for all choices of $m\in {\bf N}$ and $\{x_j\}_{j\leq m}$ in $X_K$. We refer to [MS]
and [TJ] for basic facts about type, cotype and $p$-summing operators
which are used below. The letter $c$ will always denote an absolute positive
constant, not necessarily the same in all its occurrences. By $|.|$ we
denote the cardinality of a finite set, volume of appropriate dimension,
and the Euclidean norm (this will cause no confusion). 

\section{Ellipsoidal case}

In this Section we consider the case of an arbitrary ellipsoid 
$E$ in ${\mathbf R}^n$. There exists a linear isomorphism $T:{\mathbf R}^n\rightarrow
{\mathbf R}^n$ such that $T(E)=D_n$. It will be convenient for us to write $E$
in the form
$$E= \{x=\sum_{j=1}^nx_je_j\in {\mathbf R}^n:|\sum_{j=1}^nx_ju_j|\leq 1\},\leqno (2.1)$$

\noindent where $u_j=T(e_j),\;j=1,\ldots,n$. Writing $E$ in this way, we can easily express
 $M_E$ in terms of the $u_j$'s as follows:
$$M_E=\left(\int_{S^{n-1}}\|x\|^2_{T^{-1}(D_n)}\;\sigma (dx)\right)^{1/2}
=\left(\int_{S^{n-1}}|\sum_{j=1}^nx_ju_j|^2\;\sigma (dx)\right)^{1/2}\leqno (2.2)$$
$$=\left(\frac{1}{n}\sum_{j=1}^n|u_j|^2\right)^{1/2}.$$

\medskip

Under the extra assumption that $|u_j|\leq 1,\;j=1,\ldots,n$, an
estimate for coordinate projections of $E$ was given in [Gi] in connection
with the problems of the Banach-Mazur distance to the cube and the
proportional Dvoretzky-Rogers factorization. Its
proof combines the structure of the ellipsoid with the well-known
Sauer-Shelah lemma and factorization arguments analogous to the
ones in [BT, Theorem 1.2]:

\medskip

\noindent {\bf Lemma 2.1.} {\it Let 
$E_{\tau }=\{x=\sum_{j\in\tau }x_je_j\in {\mathbf R}^{\tau }:|\sum_{j\in\tau }x_ju_j|\leq 1\}$, where $u_j\in {\mathbf R}^n$, $j\in\tau $, with $|u_j|\leq 1$.
Then, for every $\zeta\in (0,1)$ there exists $\sigma\subseteq\tau,\;|\sigma |\geq (1-\zeta )|\tau |$, such that
$${\mathcal P}_{\sigma }(E_{\tau })\supseteq c\sqrt{\zeta }\;D_{\sigma },$$
\noindent where $c>0$ is an absolute constant.$\;\;\;\Box $}

\medskip

One more step is needed in order to obtain a low $M^{\ast }$-estimate for coordinate
subspaces in the ellipsoidal case:

\medskip

\noindent {\bf Theorem 2.2.} {\it Let $E$ be an ellipsoid in ${\mathbf R}^n$. For
every $\theta\in (0,1)$ there exists a subset $\sigma $ of $\{1,\ldots,n\}$ with
$|\sigma |\geq (1-\theta )n$, such that
$${\mathcal P}_{\sigma }(E)\supseteq\frac{c\sqrt{\theta }}{\log^{1/2} (\frac{2}{\theta })M_E}
\;D_{\sigma },$$
\noindent where $c>0$ is an absolute constant.}

\medskip

\noindent {\it Proof:} We write $E$ in the form (2.1) and assume as we may that $M_E=1$. If $\rho =\{j\leq n:|u_j|\geq\sqrt{2/\theta }\}$, then by (2.2) we have $2|\rho |/\theta \leq
\sum_{j\leq n}|u_j|^2=n$,
hence $|\rho |\leq\theta n/2$. Consider the sets
of indices:

\medskip

$\tau_0=\{j\leq n:|u_j|\leq 1\}$,

$\tau_k=\{j\leq n:e^{k-1}<|u_j|\leq e^k\},\;\;\;k\geq 1$.

\medskip

\noindent If $k_0=[\log (\sqrt{2/\theta })]+1$, we have $|\bigcup_{0\leq k\leq k_0}\tau_k|\geq n-|\rho |\geq (1-\frac{\theta}{2})n$.

\medskip

We define $\zeta_k=\frac{\theta n}{2}\frac{e^k/\sqrt{|\tau_k|}}{\sum_ke^k\sqrt{|\tau_k|}}$ for all $k\leq k_0$ with $\tau_k\neq\emptyset $, and consider
the set $I=\{k\leq k_0:\tau_k\neq\emptyset\;{\rm and}\;\zeta_k<1\}$.
For each $k\in I$ we can apply Lemma 2.1 for the ellipsoid $E_{\tau_k}=
E\cap {\mathbf R}^{\tau_k}$ to find $\sigma_k\subseteq\tau_k$ with
$|\sigma_k|\geq (1-\zeta_k)|\tau_k|$ such that
$${\mathcal P}_{\sigma_k}(E_{\tau_k})\supseteq c_1\frac{\sqrt{\zeta_k}}{e^k}\;D_{\sigma_k},\leqno (2.3)$$
where $c_1$ is the constant from Lemma 2.1. Finally, we set $\sigma =\bigcup_{k\in I}\sigma_k$. Note that the above
choice of $\zeta_k$'s implies that
$$|\bigcup_{k=0}^{k_0}\tau_k|-|\sigma |\;\;\leq\;\;\sum_{k=0}^{k_0}\zeta_k|\tau_k|\;\;=\;\;
\frac{\theta n}{2},$$
and therefore, $|\sigma |\geq (1-\theta )n$.

\medskip

Suppose that $w\in D_{\sigma }$. If we write $w=\sum_{k\in I}w_k$,
where $w_k={\mathcal P}_{\sigma_k}(w)$, then by (2.3),
$$w\;\;\in\;\;\frac{1}{c_1}\sum_{k\in I}|w_k|\frac{e^k}{\sqrt{\zeta_k}}\;{\mathcal P}_{\sigma_k}(E\cap {\mathbf R}^{\tau_k})\;\;\subseteq\;\;\frac{1}{c_1}\left(\sum_{k\in I}|w_k|
\frac{e^k}{\sqrt{\zeta_k}}\right){\mathcal P}_{\sigma }(E),\leqno (2.4)$$
and since $w\in D_{\sigma }$ was arbitrary, an application of the
Cauchy-Schwartz inequality shows that
$$D_{\sigma }\subseteq\frac{1}{c_1}\left(\sum_{k\in I}\frac{e^{2k}}{\zeta_k}
\right)^{1/2}{\mathcal P}_{\sigma }(E).\leqno (2.5)$$

Inserting our $\zeta_k$'s in the sum above, we conclude that
$$D_{\sigma }\subseteq\frac{1}{c_2\sqrt{\theta n}}\left(\sum_{k=0}^{k_0}
e^k\sqrt{|\tau_k|}\right){\mathcal P}_{\sigma }(E).\leqno (2.6)$$

\medskip

It remains to give an upper bound for the sum $\sum_{k\leq k_0}e^k\sqrt{|\tau_k|}$: to this end, note that for $k=1,\ldots,k_0$, we have $|\tau_k|e^{2k-2}\leq\sum_{j\in\tau_k}|u_j|^2\leq n$
and thus $e^k\sqrt{|\tau_k|}\leq e\sqrt{n}$ for $k=1,\ldots,k_0$ which allows a
first upper bound of the order of $k_0\sqrt{n}$.

\medskip

\noindent We partition the set of indices $\{0,1,\ldots,k_0\}$ setting

\medskip

$\varphi_0=\{k\leq k_0:|\tau_k|\leq\frac{1}{k_0}\frac{n}{e^{2k-2}}\}$,

$\varphi_s=\{k\leq k_0:\frac{e^{s-1}}{k_0}\frac{n}{e^{2k-2}}<
|\tau_k|\leq\frac{e^s}{k_0}\frac{n}{e^{2k-2}}\},\;\;s\geq 1.$

\medskip

\noindent If $s_0=[\log k_0]+2$, we have $\bigcup_{0\leq s\leq s_0}\varphi_s=
\{0,1,\ldots,k_0\}$, and for every $s=1,\ldots,s_0$ we easily check that
$$|\varphi_s|\frac{e^{s-1}}{k_0}\frac{n}{e^{2k-2}}e^{2k-2}\;\;\leq\;\;\sum_{k\in\varphi_s}
\sum_{j\in\tau_k}|u_j|^2\;\;\leq \;\;n,$$
which means that
$$|\varphi_s|\leq\frac{k_0}{e^{s-1}},\leqno (2.7)$$
for all $s\leq s_0$. By the definition of $\varphi_s$ and by
(2.7), we can now estimate the sum in (2.6) as follows: 
$$\sum_{k=0}^{k_0}e^k\sqrt{|\tau_k|}\;\;=\;\;\sum_{s=0}^{s_0}\sum_{k\in\varphi_s}
e^k\sqrt{|\tau_k|}\;\;\leq\;\;\sum_{s=0}^{s_0}|\varphi_s|
\frac{e^ke^{s/2}\sqrt{n}}{\sqrt{k_0}e^{k-1}}\leqno (2.8)$$
$$\leq\;\;\frac{e\sqrt{n}}{\sqrt{k_0}}\sum_{s=0}^{s_0}\frac{k_0}{e^{s-1}}e^{s/2}
\;\;\leq \;\;e^2(\sum_{s=0}^{\infty }e^{-s/2})\sqrt{n}\sqrt{k_0}\;\;\leq \;\;c_3\sqrt{k_0}\sqrt{n}.$$
Therefore, (2.6) becomes 
$$D_{\sigma }\subseteq\frac{1}{c_4\sqrt{\theta }}\sqrt{k_0}{\mathcal P}_{\sigma } (E),\leqno (2.9)$$
which completes the proof, since $k_0\simeq\log (2/\theta )$ and we
had assumed that $M_E=1$.    $\Box $

\bigskip

We proceed to prove an extension of Theorem 2.2 concerning the case where
$E$ is an ellipsoid of dimension $m<n$ living in an arbitrary $m$-dimensional
subspace $F$ of ${\mathbf R}^n$. If $m$ is proportional to $n$, with $m/n$ sufficiently close to 1, then we still have
coordinate projections of $E$ of large dimension containing large Euclidean
balls. This result will be useful for our treatment of the general case
in Sections 3 and 4:  

\medskip

\noindent {\bf Theorem 2.3.} {\it 
Let $\varepsilon\in (0,1)$ and $F$ be a subspace of ${\mathbf R}^n$ with $\dim F=m\geq (1-\varepsilon )n$. Then, for every non-degenerate ellipsoid $E$ in $F$ and for every $\zeta\in [c_1\varepsilon\log (\frac{2}{\varepsilon }),1)$ there exists $\sigma\subseteq\{1,\ldots,n\}$ with
$|\sigma |\geq (1-\zeta )n$, such that
$${\mathcal P}_{\sigma }(E)\supseteq\frac{c\sqrt{\zeta }}{2\sqrt{2}\log^{1/2}(\frac{2}{\zeta })M_E}\;D_{\sigma },$$
\noindent where $c$ is the constant from Theorem 2.2 and $c_1=\max \{\frac{8}{c^2},\frac{1}{\log 2}\}$.}

\medskip

\noindent {\it Proof:} 
Suppose that an ellipsoid $E$ is given in $F$. We can find an orthonormal
basis $\{w_1,\ldots,w_m\}$ of $F$
and $\lambda_1,\ldots,\lambda_m>0$ such that
$$E=\{x\in F:\sum_{j=1}^m\frac{\langle x,w_j\rangle^2}{\lambda_j^2}\leq 1\}.$$
We extend to an orthonormal basis $\{w_j\}_{j\leq n}$ of
${\mathbf R}^n$ and consider the ellipsoid
$$E'=\{x\in {\mathbf R}^n:\sum_{j=1}^m\frac{\langle x,w_j\rangle^2}{\lambda_j^2}
+\sum_{j=m+1}^n\frac{\langle x,w_j\rangle^2}{b^2}\leq 1\},$$
where $b=\sqrt{\varepsilon }/M_E$. It is easy to check that
$$M_{E'}^2=\frac{1}{n}\left[\sum_{j=1}^m\frac{1}{\lambda_j^2}+\frac{n-m}{b^2}
\right]=\frac{mM_E^2+(n-m)M_E^2/\varepsilon }{n}\leq 2M_E^2.\leqno (2.10)$$

Let $\zeta\in [c_1\varepsilon\log(\frac{2}{\varepsilon }),1)$. Applying Theorem 2.2 for $E'$ and
taking into account (2.10), we find $\sigma\subseteq\{1,\ldots,n\}$ with
$|\sigma |\geq (1-\zeta)n$ for which
$${\mathcal P}_{\sigma }(E')\supseteq\frac{c\sqrt{\zeta }}{\sqrt{2}\log^{1/2} (2/\zeta)M_E}\;
D_{\sigma }.\leqno (2.11)$$

Since $\zeta\geq c_1\varepsilon\log (\frac{2}{\varepsilon})$ and the
function $\zeta /\log (\frac{2}{\zeta })$ is increasing on (0,1), one
can easily check that
$$\frac{c\sqrt{\zeta }}{\sqrt{2}\log^{1/2}(\frac{2}{\zeta })}\geq
2\sqrt{\varepsilon }.\leqno (2.12)$$
 
On the other hand, we clearly have $E'\subseteq E+bD_n$ and hence
${\mathcal P}_{\sigma }(E')\subseteq {\mathcal P}_{\sigma }(E)+bD_{\sigma }$.
Combining this with (2.11) and (2.12) we conclude that
$$\frac{c\sqrt{\zeta }}{\sqrt{2}\log^{1/2}(\frac{2}{\zeta })M_E}\;D_{\sigma }\subseteq {\mathcal P}_{\sigma }(E)+\frac{1}{2}\frac{c\sqrt{\zeta }}{\sqrt{2}\log^{1/2}(\frac{2}{\zeta })M_E}\;D_{\sigma }.\leqno (2.13)$$

\medskip

\noindent {\it Claim:} If $A$ and $B$ are convex symmetric bodies in ${\mathbf R}^{\sigma }$ and $A\subseteq B+\frac{1}{2}A$, then $A\subseteq 2B$.

[One easily checks that $A\subseteq (1+\frac{1}{2}+\ldots +\frac{1}{2^k})B+
\frac{1}{2^k}A$ and the claim follows by letting $k\rightarrow\infty $.]

\medskip

\noindent Our claim and (2.13) imply that 
$${\mathcal P}_{\sigma }(E)\supseteq\frac{c}{2\sqrt{2}}\frac{\sqrt{\zeta }}{\log^{1/2}(\frac{2}{\zeta })M_E}\;D_{\sigma },$$
and the proof of the theorem is complete.    $\Box $

\bigskip

Our next result concerns coordinate sections of ellipsoids: again, we
are interested in finding large balls contained in them. Using a result of [AM] which was recently improved in [T] (in our case each of them works
equally well), we can give an
essentially optimal answer to this question when the dimension of the
coordinate sections is small (of order roughly not exceeding $\sqrt{n}$):

\medskip

\noindent {\bf Theorem 2.4.} {\it Let $E$ be an ellipsoid in ${\mathbf R}^n$.
For every $m\leq c\sqrt{n}$ we can find a subset $\sigma $ of $\{1,\ldots,n\}$
of cardinality $|\sigma |=m$, such that
$$E\cap {\mathbf R}^{\sigma }\supseteq\frac{c'}{\sqrt{m}M_E}\;D_{\sigma }.$$
In the statement above, $c$ and $c'$ are absolute positive
constants.}

\medskip

\noindent {\it Proof:} We write $E$ in the form (2.1). As a consequence
of (2.2), observe that for every $s\leq n$ the following identity holds:
$$\Ave _{|\tau |=s}M^2_{E\cap {\mathbf R}^{\tau }}=
\left[\binom{n-1}{ s-1}/\binom{n}{s}\right]\frac{1}{s}\sum_{j=1}^n|u_j|^2
=M_E^2,\leqno (2.14)$$
where the average is over all $\tau\subseteq\{1,\ldots,n\}$ with
$|\tau |=s$. This means in particular that for every $s\leq n$ we can find
$\tau $ with $|\tau |=s$ for which $M_{E\cap {\mathbf R}^{\tau }}\leq M_E$.

\medskip

Assume that $m\leq c\sqrt{n}$ is given, where $c>0$ is an absolute constant
to be chosen. We choose $s=[\frac{m^2}{c^2}]$ and find $\tau $ with
$|\tau |=s$ and $M_{E\cap {\mathbf R}^{\tau }}\leq M_E$. Observe that
$$\Ave _{\varepsilon_j=\pm 1}\|\sum_{j\in\varphi }\varepsilon_je_j\|_E\leq
\sqrt{|\tau |}M_{E\cap{\mathbf R}^{\tau }}\leq \sqrt{|\tau |}M_E.$$

Hence, if $c$ is small
enough, the results of [AM] or [T] allow us to find $\varphi\subseteq\tau $
with $|\varphi |=2m$ such that
$$\|\sum_{j\in\varphi }\varepsilon_je_j\|_E\leq 
c_1\sqrt{|\tau |}M_E,\leqno (2.15)$$
for all $(\varepsilon_j)_{j\in\varphi }\in\{-1,1\}^{\varphi }$, where
$c_1$ is a positive absolute constant. In
other words, the coordinate section of $E$ by ${\mathbf R}^{\varphi }$ satisfies
$$E\cap {\mathbf R}^{\varphi }\supseteq\frac{1}{c_1\sqrt{|\tau |}M_E}\;Q_{\varphi }.
\leqno (2.16)$$

This means that the
identity operator $id:\ell_{\infty }^{\varphi }\rightarrow X_E\cap {\mathbf R}^{\varphi }$ has norm $\|id\|\leq c_1\sqrt{|\tau |}M_E$, and this implies that $\pi_2(id)\leq c_1K_G
\sqrt{|\tau |}M_E$ where $K_G$ is Grothendieck's constant. Applying
Pietch's factorization theorem we can find $(\lambda_i)_{i\in\varphi },\;
\sum_{i\in\varphi }\lambda_i^2=1$:
$$\|\sum_{i\in\varphi }t_ie_i\|_E
\leq c_1K_G\sqrt{|\tau |}M_E\;\left(\sum_{i\in\varphi }\lambda_i^2t_i^2\right)^{1/2} \leqno (2.17)$$
for every choice of reals $(t_i)_{i\in\varphi }$. By Markov's
inequality, we find $\sigma_1\subseteq\varphi $, $|\sigma_1|\geq
|\varphi |/2\geq m$, such that $|\lambda_i|\leq\frac{\sqrt{2}}{\sqrt{|\varphi |}}$
for all $i\in\sigma_1$. Then, for any $(t_i)_{i\in\sigma_1}$ we have
$$\|\sum_{i\in\sigma_1}t_ie_i\|_E\leq
c_1K_G\sqrt{|\tau |}M_E\frac{\sqrt{2}}{\sqrt{|\varphi |}}
\;\left(\sum_{i\in\sigma_1}t_i^2\right)^{1/2}.\leqno (2.18)$$
The choice of $|\tau |$ and $|\varphi |$ shows that
$$E\cap {\mathbf R}^{\sigma_1}\supseteq\frac{c'}{\sqrt{m}M_E}\;D_{\sigma_1},\leqno
(2.19)$$
for some absolute constant $c'>0$, and we conclude the proof by choosing any $\sigma\subseteq\sigma_1$
of cardinality $|\sigma |=m$.   $\Box $

\bigskip

\noindent {\bf Remark 2.5.} An iteration of the argument above shows that
one can extend the range of $m$'s for which Theorem 2.4 holds to e.g the set
$\{1,\ldots,[\sqrt{n}]\}$, with some loss in the constant $c'$. The dependence
on $m$ is sharp as it can be seen by the following example: consider the
ellipsoid
$E=\{(t_j)_{j\leq n}\in {\mathbf R}^n:|\sum t_ju_j|_{n+1}\leq 1\},$
where $u_j=e_j+e_{n+1},\;j=1,\ldots,n$, and $\{e_j\}_{j\leq n+1}$ is
the standard orthonormal basis in ${\mathbf R}^{n+1}$. Given any $\sigma\subseteq
\{1,\ldots,n\}$ with $|\sigma |=m$, we have that $(\frac{t}{\sqrt{m}},
\ldots,\frac{t}{\sqrt{m}})\in E\cap {\mathbf R}^{\sigma }$ precisely when
$(1+m)t^2\leq 1$. In particular, we must have $|t|\leq\frac{1}{\sqrt{m}}$.
This means that the largest ball contained in $E\cap {\mathbf R}^{\sigma }$
cannot have radius larger than $\frac{1}{\sqrt{m}}$. On the other hand,
observe that $M_E=\sqrt{2}$.

\medskip

The same example shows that the estimate in Theorem 2.2 is best possible
apart from the $\log^{1/2}(\frac{2}{\theta })$ term. By Lemma 2.1, this
logarithmic term can be removed if all the $u_j$'s are of about the same
Euclidean norm.

\section{General case: estimate in terms of the cotype-2 constant}

In this Section we study the general case, that is $K$ is an arbitrary symmetric convex body in ${\mathbf R}^n$, and
$\{e_j\}_{j\leq n}$ is a fixed orthonormal basis. 
We shall make use of the
maximal volume ellipsoid $E$ of $K$ and of the better information we have for
coordinate projections of ellipsoids. For this purpose we will also need
an estimate for the parameters $A_m(K)=\sup\{(|K\cap F|/|E\cap F|)^{1/m}:
\dim F=m\},\;m=1,\ldots,n$.

\medskip

It was proved in [BM] that the volume ratio $vr(K)=(|K|/|E|)^{1/n}$ of $K$
is bounded by $f(C_K)=cC_K[\log C_K]^4$, with the power of
$\log C_K$ improved to 1 in [MiP]. A third proof of the same fact
is given in [M4], where it is also shown that $vr(K)\leq ch(C_K)$, where
$h(y)=y\log 2y$, $y\geq 1$, and $c>0$ is an absolute constant. Our
first lemma is a modification of the argument presented in [M4] which
provides an estimate for $A_m(K),m\leq n$, in terms of $C_K$:

\medskip

\noindent {\bf Lemma 3.1.} {\it Let $K$ be a symmetric convex body in
${\mathbf R}^n$, and $E$ be the maximal volume ellipsoid of $K$. If $F$ is
an $m$-dimensional subspace of ${\mathbf R}^n$, then
$$\left(\frac{|K\cap F|}{|E\cap F|}\right)^{1/m}\leq ch(\sqrt{n/m}C_K),$$
where $h(y)=y\log 2y$, $y\geq 1$.}

\medskip

\noindent {\it Proof:} We may clearly assume that $E=D_n$. The proof
will be based on an iteration schema, analogous to the one in [M4].

We set $K_0=K$, $\alpha_0=n$, $\beta_0=n$, and for $j=1,2,\ldots$ we define:

\medskip

\noindent (i) $\alpha_j=\log\alpha_{j-1}=\log^{(j)}n$, the $j$-iterated
logarithm of $n$,

\noindent (ii) $\beta_j=\alpha_jM_{(K_{j-1}\cap F)^o}$,

\noindent (iii) $K_j=K\cap\beta_jD_n$.

\medskip

Note that for every $j$ the maximal volume ellipsoid of $K_j$ is $D_n$.
Also, $C_{K_j}\leq 2C_K$ and $d(X_{K_j},\ell_2^n)\leq\beta_j$.
By Sudakov's inequality [Su], [P1] the covering number of $K_{j-1}\cap F$ by
$\beta_jD_n\cap F$ can be estimated as follows:
$$N(K_{j-1}\cap F,\beta_jD_n\cap F)=N\leq\exp (c_1mM^2_{(K_{j-1}\cap F)^o}/\beta^2_j)=\exp (c_1m/\alpha^2_j),$$
and since, by Brunn's theorem, $|K_{j-1}\cap (x+\beta_jD_n\cap F)|\leq |K_{j-1}\cap\beta_jD_n\cap F|,\;x\in F,$ we have
$|K_{j-1}\cap F|\leq N|K_j\cap F|$ and hence
$$|K_{j-1}\cap F|^{1/m}\leq\exp (\frac{c_1}{\alpha_j^2})\;|K_j\cap F|^{1/m}
\leqno (3.1)$$

By well-known results of [DMT], [MP], and [P2] we have the string of inequalities
$$M_{(K_j\cap F)^o}\leq c_2\sqrt{\frac{n}{m}}M_{K^o_j}\leq c_3\sqrt{\frac{n}{m}}
T_2(X_{K_j^o})\leq c_4\sqrt{\frac{n}{m}}C_{K_j}\log (2d(X_{K_j},\ell_2^n))$$
and therefore
$$M_{(K_j\cap F)^o}\leq 2c_4\sqrt{\frac{n}{m}}C_K\log (2\beta_j).$$
It follows that the sequence $\{\beta_j\}_{j\geq 0}$ satisfies the
relation
$$\beta_{j+1}\leq 2c_4\sqrt{\frac{n}{m}}C_K\alpha_j\log (2\beta_j).\leqno (3.2)$$

We stop this procedure at the smallest $t$ for which $\alpha_t<6c_4$. Induction and (3.2) show that
$$\beta_t\leq 36c_4^2\sqrt{\frac{n}{m}}C_K\left[\log (\sqrt{\frac{n}{m}}C_K)+6c_4\right]\leq c'h(\sqrt{n/m}C_K).\leqno (3.3)$$

By (3.1) we see that
$$|K\cap F|^{1/m}\leq |K_t\cap F|^{1/m}\exp (c_1[\frac{1}{\alpha_1^2}+\ldots+\frac{1}{\alpha_t^2}])\leq c''|K_t\cap F|^{1/m},\leqno (3.4)$$
since $\sum\frac{1}{\alpha_j^2}$ is easily seen to be uniformly bounded. Taking into
account (3.3), (3.4) and the Blaschke-Santal\'{o} inequality we conclude that
$$\left(\frac{|K\cap F|}{|D_n\cap F|}\right)^{1/m}\leq
c''\left(\frac{|D_n\cap F|}{|(K_t\cap F)^o|}\right)^{1/m}\leq
c''M_{(K_t\cap F)^o}\leqno (3.5)$$
$$\leq 2c_4c''\sqrt{\frac{n}{m}}C_K\log (2c'h(\sqrt{n/m}C_K))\leq
ch(\sqrt{n/m}C_K).\;\;\;\;\Box$$

\medskip

Simple examples (see Remark 3.3) show that one cannot compare $M_K$ and $M_E$
even if $C_K$ is small: the only estimate one can give is that $M_E\leq\sqrt{n}M_K$, which is a direct consequence of the fact that $K\subseteq \sqrt{n}E$ by John's theorem. However, there exist subspaces $F$ of ${\mathbf R}^n$
of proportional dimension on which we can compare $M_K$ with $M_{E\cap F}$
reasonably well:

\medskip

\noindent {\bf Lemma 3.2.} {\it Let $E$ be the maximal volume ellipsoid of $K$.
For every $\varepsilon\in (0,1)$ there exists a subspace $F$ of ${\mathbf R}^n$ with $\dim F=m\geq (1-\varepsilon )n$ such that
$$M_{E\cap F}\leq\frac{ch(C_K)\log (\frac{2}{\varepsilon})}{\sqrt{\varepsilon }}\;M_K,$$
where $h(y)=y\log 2y,\;y\geq 1$, and $c>0$ is an absolute constant.}

\medskip

\noindent {\it Proof:} Let $\{w_1,\ldots,w_n\}$ be an orthonormal basis of
${\mathbf R}^n$ and $\lambda_1\geq\ldots\geq\lambda_n>0$ such that
$$E=\{x\in {\mathbf R}^n:\sum_{j=1}^n\frac{\langle x,w_j\rangle^2}{\lambda_j^2}\leq 1\}.$$

\noindent For $k=1,\ldots,n$, set $W_k={\rm span}\{w_k,\ldots,w_n\}$. By
Lemma 3.1 we have
$$\left(\frac{|K\cap W_k|}{|E\cap W_k|}\right)^{\frac{1}{n-k+1}}\leq c_1h(\sqrt{\frac{n}{n-k+1}}C_K).\leqno (3.6)$$
Note that $E\cap W_k\subseteq\lambda_k(D_n\cap W_k)$, and
hence
$$\left(\frac{|K\cap W_k|}{|E\cap W_k|}\right)^{\frac{1}{n-k+1}}\geq
\frac{1}{\lambda_k}\left(\frac{|K\cap W_k|}{|D_n\cap W_k|}\right)^{\frac{1}{n-k+1}}\leqno (3.7)$$
$$\geq\frac{1}{\lambda_kM_{K\cap W_k}}\geq\frac{1}{c_2\lambda_k\sqrt{\frac{n}{n-k+1}}M_K}.$$
Combining (3.6), (3.7) we obtain
$$\frac{1}{\lambda_k}\leq c_1c_2\sqrt{\frac{n}{n-k+1}}h(\sqrt{\frac{n}{n-k+1}}C_K)M_K\;\;\;,\;\;\;k=1,\ldots,n.\leqno (3.8)$$

Given $\varepsilon\in (0,1)$, let $m=[(1-\varepsilon )n]$ and set $F_m={\rm span}\{w_1,\ldots,w_m\}$. By (3.8) we can estimate
$M_{E\cap F_m}$ as follows:
$$M_{E\cap F_m}=\left(\frac{1}{m}\sum_{k=1}^m\frac{1}{\lambda_k^2}\right)^{1/2}\leqno (3.9)$$
$$\leq c_1c_2C_KM_K\left[\sum_{k=1}^m\frac{n^2}{m(n-k+1)^2}\log^2(2\sqrt{\frac{n}{n-k+1}}C_K)\right]^{1/2}$$
$$\leq c_1c_2C_K\log (2\sqrt{\frac{n}{n-m}}C_K)\sqrt{\frac{n}{n-m}}M_K\leq \frac{ch(C_K)\log (\frac{2}{\varepsilon })}{\sqrt{\varepsilon }}M_K.\;\;\;\Box $$

\bigskip

\noindent {\bf Remark 3.3.} The estimate (3.9) is essentially sharp, even if
$C_K$ is small: to see this, consider the class of bodies 
$K=K(a,b;s)=\{x\in {\mathbf R}^n:
\sum_{j\leq s}\frac{|x_j|}{a}+\sum_{j>s}\frac{|x_j|}{b}\leq 1\}$,
where $a,b$ are positive parameters and $s\in\{0,1,\ldots,n-1\}$. It
is clear that the ellipsoid of maximal volume in $K$ is $E=E(a,b;s)=\{x\in {\mathbf R}^n:\sum_{j\leq s}\frac{|x_j|^2}{a^2}+\sum_{j>s}\frac{|x_j|^2}{b^2}\leq 1\}$.
It is also clear that both the cotype-2 constant and the volume ratio
of $K$ are uniformly bounded (independently of $n,s,a$ and $b$).

\medskip

Given $\varepsilon\in (0,1)$, choose $b=a\sqrt{\varepsilon}$, $s=m=(1-\varepsilon )n$. Then, it is easy to check that 
$M_K\simeq \sqrt{n}\sqrt{\varepsilon }/a$, while
$M_{E\cap F_m}\simeq\sqrt{n}/a$.

\medskip

Also, we can have the ratio $M_E/M_K$ as close to $\sqrt{n}$ as we
like: choose, for example, $s=n-1$ and $b=\frac{a}{n-1}$. Then, $M_K\simeq 1/\sqrt{n}b$ while $M_E\simeq 1/b$.

\bigskip

\medskip

Combining Theorem 2.3 and Lemma 3.2 we prove our $M^{\ast}$-estimate in
terms of the cotype-2 constant of $X_K$:

\medskip

\noindent {\bf Theorem 3.4.} {\it Let $K$ be a symmetric convex body in
${\mathbf R}^n$, and $X_K=({\mathbf R}^n,\|.\|_K)$. For every $\theta\in (0,1)$ there
exists $\sigma\subseteq\{1,\ldots,n\},\;|\sigma |\geq (1-\theta )n$, such
that
$${\mathcal P}_{\sigma }(K)\supseteq\frac{c\theta }{\log^2(\frac{2}{\theta })h(C_K)M_K}\;D_{\sigma },$$
where $h(y)=y\log 2y,\;y\geq 1$, and $c>0$ is an absolute constant.}

\medskip

\noindent {\it Proof:} Let $E$ be the maximal volume ellipsoid of $K$, and set
$\varepsilon =\varepsilon (\theta )=\theta /c_2\log (\frac{2}{\theta })$,
where $c_2>0$ is a constant to be chosen. By Lemma 3.2 we can find a
subspace $F$ of ${\mathbf R}^n$ with $\dim F\geq (1-\varepsilon )n$ such
that
$$M_{E\cap F}\leq\frac{c_3h(C_K)\log (\frac{2}{\varepsilon })}{\sqrt{\varepsilon }}\;M_K.\leqno (3.10)$$

Observe that if $c_2$ is large enough, then $\theta\geq c_1\varepsilon\log(\frac{2}{\varepsilon })$ where $c_1$ is the constant in Theorem 2.3. Thus, we can
apply Theorem 2.3 for $E\cap F$ to find $\sigma\subseteq\{1,\ldots,n\}$ with
$|\sigma |\geq (1-\theta )n$ for which
$${\mathcal P}_{\sigma }(E\cap F)\supseteq\frac{c\sqrt{\theta }}{2\sqrt{2}\log^{1/2}(\frac{2}{\theta })M_{E\cap F}}\;
D_{\sigma }.\leqno (3.11)$$

Combining (3.10) with (3.11) we finish the proof.   $\Box $

\bigskip

\noindent {\bf Remark 3.5.} It should be noted that the estimate 
given by Theorem 3.4 is exact
not only when $C_K$ is small (like e.g in the ellipsoidal case),
but in the whole range $[1,\sqrt{n}]$ of possible values of $C_K$
i.e even if $C_K$ is extremely large. This can be easily seen if
one considers the case of $B_p^n,p>2$, the unit ball of $\ell_p^n$,
and the standard coordinate system in ${\mathbf R}^n$.
Fix for example $\theta =\frac{1}{2}$. Then, the radius of the
largest Euclidean ball inscribed in any $[\frac{n}{2}]$-dimensional
coordinate projection of $B_p^n$ is 1, and the well-known
estimates for $C_{B_p^n}$ and $M_{B_p^n}$ show that Theorem 3.4
is sharp apart from logarithmic terms. We do not know if the
``almost linear" dependence on $\theta $ which our method provides
is optimal. However, the ellipsoidal case shows that $\sqrt{\theta }$
dependence is the best one might hope for.

\bigskip

\noindent {\bf Remark 3.6.} One can give a weaker estimate, analogous
to the one obtained in Theorem 3.4, using the isomorphic Sauer-Shelah
lemma of Szarek-Talagrand [ST] and a factorization result of Maurey
[Ma] (see also [TJ]). Starting with the body $K$ and the orthonormal
basis $\{e_j\}_{j\leq n}$, we have the inequality 
$$\Ave _{\varepsilon_j
=\pm 1}\|\sum_{j=1}^n\varepsilon_je_j\|_K\leq\sqrt{n}M_K,$$
\noindent and therefore, by Markov's inequality we can find
${\mathcal A}\subseteq\{-1,1\}^n$ of cardinality $|{\mathcal A}|\geq 2^{n-1}$
such that $\|\sum \varepsilon_je_j\|_K\leq 2M_K\sqrt{n}$ whenever
$\varepsilon =(\varepsilon_1,\ldots,\varepsilon_n)\in {\mathcal A}$.
If we view ${\mathcal A}$ as a set of points in ${\mathbf R}^n$, this
means that ${\mathcal A}\subseteq 2M_K\sqrt{n}K$. On the other hand,
the isomorphic Sauer-Shelah lemma shows that for some absolute
constant $c_1>0$ and for every $\theta\in (0,1)$ there exists
$\sigma\subseteq\{1,\ldots,n\}$, $|\sigma |\geq (1-\frac{\theta }{2})n$,
with ${\rm co}({\mathcal P}_{\sigma }({\mathcal A})\supseteq c_1\frac{\theta }{2}Q_{\sigma }$, and hence 
$${\mathcal P}_{\sigma }(K)\supseteq \frac{c_1\theta }{4M_K\sqrt{n}}\;Q_{\sigma }.$$
It follows that if $Y=({\mathbf R}^{\sigma_1},\|.\|_{K^o})$, then
$id:\ell_{\infty }^{\sigma_1}\rightarrow Y^{\ast }$ has norm
$\|id\|\leq\frac{4M_K\sqrt{n}}{c_1\theta }$, and Maurey's theorem
shows that
$$\pi_2(id)\leq c_2\frac{M_K\sqrt{n}}{\theta }\;g(Y^{\ast }),$$
\noindent where $g(Y^{\ast })=C_{Y^{\ast }}\sqrt{1+\log (C_{Y^{\ast }})}$.
Then, we can apply Pietch's factorization theorem in the context
of [BT, Theorem 1.2] to find $\sigma\subseteq\sigma_1$ with
$|\sigma |\geq (1-\frac{\theta }{2})|\sigma_1|\geq (1-\theta )n$
for which
$$\left(\sum_{i\in\sigma }t_i^2\right)^{1/2}\leq
c_3\frac{M_Kg(Y^{\ast })}{\theta^{3/2}}\/\|\sum_{i\in\sigma }t_ie_i\|_{K^o}$$
\noindent is true for all $(t_i)_{i\in\sigma }$. Taking polars in
${\mathbf R}^{\sigma }$ and using the fact that
$C_{Y^{\ast }}\leq c_4C_K\|{\rm Rad}_{X_K}\|$, we conclude that
$${\mathcal P}_{\sigma }(K)\supseteq\frac{c\theta^{3/2}}{f(K)M_K}\;D_{\sigma },$$
\noindent where $c>0$ is an absolute constant, and
$f(K)=C_K\|{\rm Rad}_{X_K}\|\sqrt{1+\log (C_K\|{\rm Rad}_{X_K}\|)}$.

\bigskip

\noindent {\bf Remark 3.7.} One can modify the proof of Theorem 3.6
to give analogous estimates in which $M_K$ is replaced by other
``volumic'' parameters of $K$ or $K^o$.

Consider e.g the sequence of volume numbers of $K^o$
$$v_s(K^o)=\max\{(|{\mathcal P}_F(K^o)|/|D_n\cap F|)^{1/s}:\;\dim F=s\},
\leqno (3.12)$$
\noindent where $s=1,\ldots,n$. As a consequence of the
Aleksandrov-Fenchel inequalities, one can easily see that
$\{v_s(K^o)\}_{s\leq n}$ is non increasing (see [P1]):
$$v_1(K^o)\geq v_2(K^o)\geq \ldots \geq v_n(K^o)=v.rad(K^o).\leqno (3.13)$$
Let $K$ be a symmetric convex body in ${\mathbf R}^n$ and let $E$ be the
ellipsoid of maximal volume in $K$ as in Lemma 3.4. Using the inverse
Santal\'{o} inequality in (3.6), (3.7) we get
$$\frac{1}{\lambda_{k}}\leq \left (\frac{|D_n\cap W_k|}{|K\cap W_k|}
\right )^{\frac{1}{n-k+1}}\left (\frac{|K\cap W_k|}{|E\cap W_k|}\right )
^{\frac{1}{n-k+1}}\leqno (3.14)$$
$$\leq c\left (\frac{|{\mathcal P}_{W_k}(K^o)|}{|D_n\cap W_k|}\right )^{\frac{1}{n-k+1}}c_1h(\sqrt{\frac{n}{n-k+1}}C_K)$$
\noindent for $k=1,\ldots,n$. By the definition (3.11) of $v_{n-k+1}(K^o)$
this means that
$$\frac{1}{\lambda_{k}}\leq c_2h(C_K)v_{n-k+1}(K^o)\sqrt{\frac{n}{n-k+1}}
\log (2\sqrt{\frac{n}{n-k+1}}).\leqno (3.15)$$
Inserting this estimate in (3.9) we obtain:
$$M_{E\cap F_m}^2=\frac{1}{m}\sum_{k=1}^{m}\frac{1}{\lambda_{k}^2}
\leq\frac{c_2^2h^2(C_K)\log^2(\frac{2n}{n-m+1})}{m}\sum_{k=1}^{m}
\frac{n}{n-k+1}v_{n-k+1}^2(K^o).\leqno (3.16)$$
\noindent The monotonicity of volume numbers shows that
$v_{n-k+1}(K^o)\leq v_{n-m+1}(K^o),\;k=1,\ldots,m$, and combining
with the fact that
$$\sum_{k=1}^{m}\frac{n}{n-k+1}\leq n\log (\frac{n}{n-m})$$
we arrive at
$$M_{E\cap F_m}\leq \frac {cn}{m}h(C_K)v_{n-m+1}(K^o)\log^{3/2}(\frac{2n}{n-m}).\leqno (3.17)$$

Set $m=[(1-\theta )n]$. Then, (3.17) can be rewritten as
$$M_{E\cap F_m}\leq c'h(C_K)v_{[\theta n]}(K^o)\log^{3/2}(\frac{2}{\theta }),$$
\noindent and, using Theorem 2.3 exactly as in the proof of Theorem 3.4,
we can find $\sigma\subseteq\{1,\ldots,n\}$ with $|\sigma |\geq (1-c_1\theta \log (\frac{2}{\theta }))n$ for which
$${\mathcal P}_{\sigma }(K)\supseteq\frac{c\sqrt{\theta }}{\log^{3/2}(\frac{2}{\theta })v_{[\theta n]}(K^o)h(C_K)}\;D_{\sigma }.\leqno (3.18)$$

A similar argument shows that for some $\sigma $ of the same cardinality
we have
$${\mathcal P}_{\sigma }(K)\supseteq\frac{c\sqrt{\theta }w_{[\theta n]}(K)}
{\log^{3/2}(\frac{2}{\theta })h(C_K)}\;D_{\sigma },\leqno (3.19)$$
\noindent where $w_s(K)=\min\{(|K\cap F|/|D_n\cap F|)^{1/s}:\;\dim F=s\},\;
s=1,\ldots,n$.

\section {General case: estimate in terms of the volume ratio}

In this Section we use the volume ratio $vr(K)$ of $K$ instead of
the cotype-2 constant of $X_K$ as a parameter for our low
$M^{\ast }$-estimate. Let $E$ be the maximal volume ellipsoid
of $K$. We start with a lemma which estimates the covering number $N(K,E)$ in
terms of the volume ratio $vr(K)=(|K|/|E|)^{1/n}$. Our proof
is based on Lemma 4.4 from [MS2], actually the argument given there leads
to a stronger estimate, but we include a simple proof of what we
need here for the sake of completeness.
Recall that $N(K,L)$ is the smallest natural number $N$ for which
there exist $x_1,\ldots,x_N\in {\mathbf R}^n$ with $K\subseteq\bigcup_{i\leq N}(x_i+L)$:

\medskip

\noindent {\bf Lemma 4.1.} {\it Let $K$ and $L$ be two symmetric
convex bodies in ${\mathbf R}^n$ such that $L\subseteq K$. Then,
$$N(K,L)\leq c^n\frac{|K|}{|L|},$$
\noindent where $c>0$ is an absolute constant.}

\medskip

\noindent {\it Proof:} Consider a set
$N$ of points in $K$ such that $\|x-x'\|_L\geq 1$ for
every $x,x'\in N,\;x\neq x'$, which has the maximal possible
cardinality. Observe that the sets $\frac{2}{3}x+\frac{L}{3},\;x\in N$
have disjoint interiors and, since $L\subseteq K$, they are all contained in $K$.
We easily deduce that
$$|N|\leq  3^n\frac{|K|}{|L|}.\leqno (4.1)$$

\noindent Finally, it is clear that $K\subseteq\bigcup_{x\in N}(x+L)$,
which completes the proof.   $\Box $

\medskip

Suppose that $K$ is any symmetric convex body in ${\mathbf R}^n$
and $E$ is the ellipsoid of maximal volume in $K$. The analogue of Lemma 3.2
in the ``volume ratio" formulation is the following:

\medskip

\noindent {\bf Lemma 4.2.} {\it Let $E$ be the maximal volume
ellipsoid of $K$. For every $\varepsilon\in (0,1)$ there exists
a subspace $F$ of ${\mathbf R}^n$ with $\dim F=m\geq (1-\varepsilon )n$,
such that
$$M_{E\cap F}\leq
[c\;vr(K)]^{1/\varepsilon }M_K,$$
\noindent where $c>0$ is an absolute constant.}

\medskip

\noindent {\it Proof:} As in the proof of Lemma 3.2, let
$$E=\{x\in {\mathbf R}^n:\;\sum_{j=1}^n\frac{\langle x,w_j\rangle^2}{\lambda_j^2}
\leq 1 \},$$
\noindent where $\{w_1,\ldots,w_n\}$ is an orthonormal basis of
${\mathbf R}^n$ and $\lambda_1\geq\ldots\geq\lambda_n>0$. Fix
$k\in\{1,\ldots,n\}$ and consider the subspace $W_k={\rm span}\{w_k,
\ldots,w_n\}$. According to Lemma 4.1, we can find $x_1,\ldots,x_N\in K$
such that $N=N(K,E)\leq [c_1vr(K)]^n$ and $K\subseteq\bigcup (x_i+E)$.
Project all the $(x_i+E)$'s onto $W_k$. Then,
$$K\cap W_k\subseteq{\mathcal P}_{W_k}(K)\subseteq
\bigcup_{j\leq N}{\mathcal P}_{W_k}(x_j+E)=
\bigcup_{j\leq N}({\mathcal P}_{W_k}(x_i)+E\cap W_k),\leqno (4.2)$$
and hence, $N(K\cap W_k, E\cap W_k)\leq N(K,E)$. Thus, we can
estimate the ratio of the volumes of $K\cap W_k$ and $E\cap W_k$ using
(4.2):
$$\left(\frac{|K\cap W_k|}{|E\cap W_k|}\right)^{\frac{1}{n-k+1}}\leq
[N(K,E)]^{\frac{1}{n-k+1}}\leq [c_1vr(K)]^{\frac{n}{n-k+1}}.\leqno (4.3)$$
Combining with (3.7) we get
$$\frac{1}{\lambda_k}\leq c_2\sqrt{\frac{n}{n-k+1}}\;[c_1vr(K)]^{\frac{n}{n-k+1}}M_K,\;\;\;\;k=1,\ldots,n.\leqno (4.4)$$
We continue as in Lemma 3.2: Given any $\varepsilon\in (0,1)$,
we consider the first $m$ for which $m\geq (1-\varepsilon )n$ and set
$F_m={\rm span}\{w_1,\ldots,w_m\}$. In view of (4.5) we can compare
$M_{E\cap F_m}$ with $M_K$ as follows:
$$M_{E\cap F_m}=\left(\frac{1}{m}\sum_{k=1}^m\frac{1}{\lambda_k^2}\right)^{\frac{1}{2}}\leqno (4.5)$$
$$\leq c_2M_K\left(\sum_{k=1}^m\frac{n}{m(n-k+1)}[c_1vr(K)]^{\frac{2n}{n-k+1}}
\right)^{\frac{1}{2}}$$
$$\leq M_K[c_3vr(K)]^{\frac{n}{n-m+1}}\sqrt{\frac{n}{m}}\log^{1/2}(\frac{n}{n-m}),$$
and the lemma follows with the observation that
$\frac{\log (1/\varepsilon )}{1-\varepsilon }\rightarrow 1$ as $\varepsilon\rightarrow 1^-$.   $\Box $

\medskip

\noindent {\bf Remark 4.3.} By well-known results of S.J. Szarek
and N. Tomczak-Jaegermann (see [Sz], [STJ]) which were extending
previous work of Kashin, if $E$ is the maximal volume ellipsoid
of $K$, then for every $k=1,\ldots,n-1$ there exist $k$-dimensional 
subspaces $F$
of ${\mathbf R}^n$ for which $E\cap F\subseteq K\cap F\subseteq
(c\;vr(K))^{\frac{n}{n-k}}E\cap F$, and this obviously implies
that $M_{E\cap F}\leq [c\;vr(K)]^{\frac{n}{n-k}}M_{K\cap F}$.
This leads to the same estimate as in Lemma 4.2 above, actually
if $E=D_n$ this is true for all subspaces $F$ in a 
subset ${\mathcal A}$ of $G_{n,k}$
with almost full measure $\nu_{n,k}({\mathcal A})>1-2^{-n}$.
The argument provided by Lemmata 4.1 and 4.2 gives a concrete
example of a subspace on which the weaker ``$M_{E\cap F}$ and
$M_{K\cap F}$" comparison is true: it can be chosen as the
$k$-dimensional subspace which is coordinate with respect
to $E$ and corresponds to the $k$ largest semiaxes of $E$.
If $E=D_n$, then this weak comparison is true for all $F\in G_{n,k}$.

\medskip

Combining Lemma 4.2 with Theorem 2.3 we prove our volume-ratio result:

\medskip

\noindent {\bf Theorem 4.4.} {\it Let $K$ be a symmetric convex body in ${\mathbf R}^n$. For every $\theta\in (0,1)$ there exists $\sigma\subseteq\{1,\ldots,n\},\;
|\sigma |\geq (1-\theta )n$, such that
$${\mathcal P}_{\sigma }(K)\supseteq\frac{1}{[c_1vr(K)]^{\frac{c_2\log(\frac{2}{\theta })}{\theta }}M_K}\;D_{\sigma },$$
\noindent where $c_1,c_2$ are absolute positive constants}.

\medskip

\noindent {\it Proof:} Let $E$ be the maximal volume ellipsoid of $K$,
and set $\varepsilon =\varepsilon (\theta )=\frac{\theta }{c_2\log (\frac{2}{\theta })}$, where $c_2>0$ is a constant to be chosen. Using Lemma 4.2 we
find a subspace $F$ of ${\mathbf R}^n$ with $\dim F\geq (1-\varepsilon )n$
such that
$$M_{E\cap F}\leq [c_4vr(K)]^{1/\varepsilon }
M_K.\leqno (4.6)$$
If $c_2$ is large enough, we easily check that
$\theta\geq c_1\varepsilon\log (\frac{2}{\varepsilon })$ where
$c_1$ is the constant in Theorem 2.3. We can therefore apply
Theorem 2.3 for $E\cap F$ to find $\sigma\subseteq\{1,\ldots,n\}$
with $|\sigma |\geq (1-\theta )n$, such that
$${\mathcal P}_{\sigma }(E\cap F)\supseteq\frac{c\sqrt{\theta }}{2\sqrt{2}
\log^{1/2}(\frac{2}{\theta })M_{E\cap F}}\;D_{\sigma }.\leqno (4.7)$$
Combining (4.6) with (4.7) we conclude the
proof.    $\Box $

\medskip

For classes of spaces with uniformly bounded volume ratio, Theorem 4.3
gives an optimal answer as long as, say, $\theta\geq\frac{1}{2}$.
The estimate obtained ``explodes" if $vr(K)$ is large or if $\theta $
is needed to be close to 0.

\section{Linear duality relations for coordinate
sections of ellipsoids}

Let $K$ be a symmetric convex body in ${\mathbf R}^n$. We introduce the
integer valued functions $t,t_c:{\mathbf R}^+\rightarrow {\bf N}$ defined by
$$t(r)=t(K,r)=\max\{k\leq n: {\rm there\;exists\;a\;subspace}\;E\;$$
$${\rm with\;dim}E=k,
{\rm such\;that}
\frac{1}{r}|x|\leq\|x\|\;{\rm for\;every}\;x\in E\}$$
and
$$t_c(r)=t_c(K,r)=\max\{k\leq n:{\rm there\;exists\;a\;coordinate\;subspace}\;E\;$$
$${\rm with}
\;\dim E=k
\;{\rm such\;that}\;\frac{1}{r}|x|\leq\|x\|\;{\rm for\;
every\;}x\in E\}.$$

It is easy to see that if $K$ is an ellipsoid in ${\mathbf R}^n$, then
$t(K,r)+t(K^o,\frac{1}{r})\geq n$. In [M5] it is proved that for
every body $K$, for every $r>0$, and for every $\tau\in (0,1)$, one
has a similar duality relation:
$$t(K,r)+t(K^o,\frac{1}{\tau r})\geq (1-\tau )n-C,\leqno (5.1)$$
\noindent where $C>0$ is a universal constant. The proof of this
fact is based on the strong form (1.2) of the low $M^{\ast }$-estimate
and on the ``distance lemma": if $\frac{1}{a}|x|\leq\|x\|\leq b|x|$
for every $x\in {\mathbf R}^n$ and if $(M_K/b)^2+(M_{K^o}/a)^2=s>1$,
then $ab\leq\frac{1}{s-1}$.

\medskip

In this Section we establish a coordinate version of (5.1) in
the ellipsoidal case. Our estimate depends on how close the ellipsoid
is to being in M-position:

\medskip

\noindent {\it Definition:} For a symmetric convex body $K$ in ${\mathbf R}^n$
we denote by $\lambda_K$ its volume radius: $\lambda_K=(|K|/|D|)^{1/n}$.
We also write $N_K$ for $N(K,\lambda_KD)$ and say that $K$ is in
$M_{\delta }$-{\it position} if $\delta\geq\frac{1}{n}\log N_K$.

\medskip 

Our first lemma provides some simple estimates
which show that this position is ``stable'' under the operations of
taking intersection or convex hull with a ball:

\medskip

\noindent {\bf Lemma 5.1.} {\it Let $K$ be a symmetric convex body
in ${\mathbf R}^n$, and let $r, r_1>0$ be given. Define $K_r=K\cap rD$
and $K^{r_1}={\rm co}(K\cup r_1D)$. Then,

{\rm (i)} $N_{K_r}\leq\max\{3^n N_K^2, 9^n N_K\}$.

{\rm (ii)} $N_{K^{r_1}}\leq 5^n N_K$.} 

\medskip

\noindent {\it Proof:} (i) From the Brunn-Minkowski inequality
it easily follows that $|K\cap rD|\geq |K\cap (x+rD)|,\;x\in {\mathbf R}^n$.
This implies that $|K|\leq N(K,rD)|K\cap rD|$ or, equivalently,
$$\lambda_K^n\leq N(K,rD)\lambda_{K_r}^n.\leqno (5.2)$$
\noindent We distinguish two cases:

(1) If $\lambda_K<r$, then $N(K,rD)\leq N_K$ and, by (5.2),
$\lambda_K^n\leq N_K\lambda_{K_r}^n$. It follows that
$$N_{K_r}\leq N(K,\lambda_{K_r}D)\leq N_KN(D,\frac{\lambda_{K_r}}{\lambda_K}D)
\leq N_K N(D,\frac{1}{N_K}D)\leq 3^n N_K^2.$$

(2) If $\lambda_K>r$, then $N(K,rD)\leq N(K,\lambda_KD)N(D,\frac{r}{\lambda_K}D)\leq N_K3^n(\frac{\lambda_K}{r})^n$ and hence, by (5.2),
$(\frac{r}{\lambda_{K_r}})^n\leq 3^nN_K$. It follows that
$$N_{K_r}\leq N(rD,\lambda_{K_r}D)\leq 3^n(\frac{r}{\lambda_{K_r}})^n
\leq 9^nN_K.$$

\medskip

\noindent (ii) We obviously have $\lambda_{K^{r_1}}\geq\max\{\lambda_K,r_1\}$.
Also, $K^{r_1}\subseteq K+r_1D$, which gives
$$N_{K^{r_1}}\leq N(K^{r_1}, 2\lambda_{K^{r_1}}D)N(D,\frac{1}{2}D)\leq
5^nN(K+r_1D,(\lambda_K+r_1)D)\leq 5^nN_K.\;\;\;\Box $$

\medskip

For an arbitrary symmetric convex body $K$, one has in general the
information $\lambda_K M_K\geq 1$ as a consequence of the polar
coordinates formula for volume. Our next lemma provides an ``inverse"
inequality in terms of the parameters $N_{K^o}$ and $b=\sup\{\|x\|: x\in S^{n-1}\}$:

\medskip

\noindent {\bf Lemma 5.2.} {\it Let $K$ be a symmetric convex body in
${\mathbf R}^n$, and assume that $\|x\|\leq b|x|$ for all $x\in {\mathbf R}^n$.
Then,
$$M_K\leq\frac{c}{\lambda_K}N_{K^o}^{t/n}$$
\noindent where $c>0$ is an absolute constant, and
$t\leq C(\frac{b}{M_K})^2$.}

\medskip

\noindent {\it Proof:} Using Theorem 6 from [BLM] (to be more
precise, using an argument identical to the one given there
and the observation that what is really used is the ratio
$b/M_K$), one can find
orthogonal transformations $u_1,\ldots,u_t\in O(n)$ such that
$$\frac{M_K}{2}D\subseteq\;T=\frac{1}{t}\sum_{i=1}^tu_i(K^o)\;
\subseteq 2M_K D,\leqno (5.4)$$
\noindent with $t\leq C(\frac{b}{M_K})^2$, where $C>0$ is an
absolute constant.

On observing that $N(T,\lambda_{K^o}D)\leq [(N(K^o,\lambda_{K^o}D)]^t
=N_{K^o}^t$, we can estimate $M_K$ by (5.4) as follows:
$$M_K\leq 2(\frac{|T|}{|D|})^{1/n}\leq 2\lambda_{K^o}N_{K^o}^{t/n}.
\leqno (5.5)$$
Finally, the Blaschke-Santal\'{o} inequality implies that
$\lambda_K\lambda_{K^o}\leq 1$, and hence the proof of the
Lemma is complete.   $\Box $

\medskip

We can now pass to the proof of the main result of this section:

\medskip

\noindent {\bf Theorem 5.3.} {\it Let $E$ be an ellipsoid in ${\mathbf R}^n$,
and assume that both $E$ and $E^o$ are in $M_{\delta }$-position.
For every $r>0$ and every $\tau\in (0,1)$ we have
$$t_c(E,r)+t_c(E^o,\frac{u(\tau ,\delta )}{r})\geq (1-\tau )n,$$
\noindent where $u(\tau ,\delta )=\frac{c\log (\frac{2}{\tau })}{\sqrt{\tau }}
e^{\frac{c\delta\log^2 (\frac{2}{\tau })}{\tau }}$, and $c>0$ is
an absolute constant.}

\medskip

\noindent {\it Proof:} Let $r>0$ and $\tau\in (0,1)$ be given. Consider
the body $E_r=E\cap rD$. Since $E_r$ is $\sqrt{2}$-isomorphic to an
ellipsoid, one can easily check that Theorem 2.2 holds for $E_r$: for
every $\theta\in (0,1)$ we can find $\sigma\subseteq\{1,\ldots,n\}$
with $|\sigma |\geq (1-\theta )n$ such that ${\mathcal P}_{\sigma }(E_r^o)
\supseteq [g(\theta )/M(E_r^o)]D_{\sigma }$, where $g(\theta )=c\sqrt{\theta }
/2\sqrt{\log (2/\theta )}$ and $c$ is the same constant as in Theorem 2.2.

\medskip

\noindent We distinguish three cases:

\medskip

\noindent {\it Case 1:} $\frac{M(E_r^o)}{r}\in [g(\tau ), g(1))$.

\medskip

In this case, consider any $\lambda\in (\tau , 1]$ with $\frac{1}{r}M(E_r^o)<
g(\lambda )$. We can find $\sigma_1\subseteq\{1,\ldots,n\}$
with $|\sigma_1|\geq (1-\lambda )n$ such that
$${\mathcal P}_{\sigma_1}(E_r^o)\supseteq\frac{g(\lambda )}{M(E_r^o)}D_{\sigma_1}
,$$
\noindent and it is easy to check that, for every $x\in {\mathbf R}^{\sigma_1}$, 
$\max\{\|x\|,\frac{1}{r}|x|\}=\|x\|_{E_r}>\frac{1}{r}|x|$, which means that
$\frac{1}{r}|x|\leq\|x\|$, i.e
$$t_c(E,r)\geq (1-\lambda )n.\leqno (5.6)$$
\noindent Taking the infimum of all $\lambda $'s for which $\frac{M(E_r^o)}{r}<g(\lambda )$, we conclude that (5.6) also holds for the solution in $\lambda $
of the equation $M(E_r^o)=rg(\lambda )$.

\medskip

Now, choose $\mu\in (0,1)$ such that $(1-\lambda )+(1-\mu )=1-\tau $,
and $r_1>0$ satisfying $M((E_r)^{r_1})r_1<g(\mu )$ (this is always possible
since the left hand side is decreasing in $r_1$ and tends to zero as
$r_1\rightarrow\infty $). Since $(E_r)^{r_1}$ is 2-isomorphic to an
ellipsoid, we can find $\sigma_2\subseteq\{1,\ldots,n\}$, $|\sigma_2|\geq
(1-\mu )n$, with
$${\mathcal P}_{\sigma_2}((E_r)^{r_1})\supseteq \frac{g(\mu )}{M((E_r)^{r_1})}D_{\sigma_2}
,$$
\noindent thus $\max\{r_1|x|,\|x\|_{E_r^o}\}=\|x\|_{[(E_r)^{r_1}]^o}\geq\frac
{g(\mu )}{M((E_r)^{r_1})}|x|>r_1|x|$, i.e $\|x\|_{E^o}\geq\|x\|_{E_r^o}>r_1|x|$ on
${\mathbf R}^{\sigma_2}$, which means that
$$t_c(E^o,\frac{1}{r_1})\geq (1-\mu )n.\leqno (5.7)$$
\noindent Again, we may take $r_1$ to be the solution of the equation
$M((E_r)^{r_1})r_1=g(\mu )$ in $r_1$.

\medskip

Combining (5.6) with (5.7) we obtain
$$t_c(E,r)+t_c(E^o,\frac{1}{r_1})\geq (1-\tau )n,\leqno (5.8)$$
\noindent and it remains to compare $r$ with $r_1$. Let us write $W$ for
the body $(E_r)^{r_1}$. By the way $W$ has been constructed, it is easily
checked that the following are satisfied:

\medskip

(i) $M(W)r_1=g(\mu )$ and $M(W^o)\geq M(E_r^o)=rg(\lambda )$.

\medskip

(ii)$\|x\|_W\leq\frac{1}{r_1}|x|$ and $\|x\|_{W^o}\leq r|x|$, $x\in {\mathbf R}^n$.

\medskip

(iii) $N_W^{1/n}\leq c_1N_E^{c_2/n}$ and $N_{W^o}^{1/n}\leq
c_1N_{E^o}^{c_2/n}$, where $c_1,c_2>0$ are absolute constants. This
is a simple consequence of Lemma 5.1, since both $W$ and $W^o$ are
formed from $E$ and $E^o$ with two successive operations of taking
intersection and convex hull with balls.

\medskip

We simply write
$$\frac{r}{r_1}=\frac{r}{M(W^o)}\;\frac{1}{r_1M(W)}\;M(W) M(W^o)$$
\noindent and making use of (i)-(iii) and of Lemma 5.2 we arrive at
$$\frac{r}{r_1}\leq \frac{c}{g(\lambda )g(\mu )}N_{E^o}^{C/ng^2(\mu )}
N_E^{C/ng^2(\lambda )}.\leqno (5.9)$$
\noindent Note that, at some point, we also used the fact that
$\lambda_E\lambda_{E^o}\simeq 1$. Finally, assuming that both $E$
and $E^o$ are in $M_{\delta }$-position, we rewrite (5.9) as follows:
$$\frac{r}{r_1}\leq\frac{c}{g(\lambda )g(\mu )}e^{C\delta /g^2(\lambda )g^2(\mu )}.\leqno (5.10)$$
\noindent We have $\lambda +\mu =1+\tau $ and with this condition
we can easily check that $\frac{1}{g(\lambda )g(\mu )}\leq\frac{c\log (\frac{2}{\tau })}{\sqrt{\tau }}$, which completes the proof in this case.

\medskip

\noindent {\it Case 2:} $\frac{M(E_r^o)}{r}\geq g(1)$.

\medskip

We choose $r_1>0$ such that $M((E_r)^{r_1})r_1=g(\tau )$ and as above we
conclude that $t_c(E^o,\frac{1}{r_1})\geq (1-\tau )n$. The estimate
for $r/r_1$ is done exactly in the same way, the only difference
being that now $r/M(E_r^o)\leq 1/g(1)$.

\medskip

\noindent {\it Case 3:} $\frac{M((E_r)^o)}{r}<g(\tau )$.

\medskip

This is the simplest case since we already have $t_c(E,r)\geq (1-\tau )n$.
  $\Box $

\section{Integer points inside an ellipsoid: some remarks}

Consider an arbitrary ellipsoid $E$ in ${\mathbf R}^n$. Write $E$ in the
form (2.1), so that $\sum_{j\leq n}|u_j|^2=nM_E^2$. Without loss of
generality we may assume that the $|u_j|$'s are arranged in the
increasing order, therefore a simple application of Markov's inequality
shows that
$$|u_j|\leq\sqrt{\frac{n}{n-j+1}}M_E\;\;\;,\;\;\;j=1,\ldots,n.\leqno (6.1)$$
Recall that the $j$-th successive minimum $\lambda_j(E)$ of $E$ is
defined by $\lambda_j(E)=\min\{\lambda >0: {\rm dim }({\rm span }(\lambda E
\cap {\bf Z}^n))\geq j\}$. Inequality (6.1) gives an estimate on the
successive minima of $E$ in terms of $M_E$:

\medskip

\noindent {\bf Fact I}: {\it Let $E$ be an ellipsoid in ${\mathbf R}^n$. Then,
$\lambda_j(E)\leq\sqrt{\frac{n}{n-j+1}}M_E,\;j=1,\ldots,n$. In particular,
if $M_E\leq 1$ then $E$ contains an integer point different from the
origin}.

\medskip

Note that if $M_E>1$ then $E$ may contain no integer points other
than the origin. Consider for example a ball of radius $r=\frac{1}{M_E}$.

\medskip

Let us concentrate on the case $M_E<1$. If $M_E<|D_n|^{1/n}/2$, then
we obviously have $|E|>2^n$ and Minkowski's theorem with its
relatives start giving estimates on the cardinality of the set of
integer points in $E$. We are interested in the range
$|D_n|^{1/n}/2<\;M_E<\;1$. From Fact I we know that $E$ contains
non-trivial integer points, and using $M_E$ as a parameter we
try to estimate the number of them. Theorem 2.4 can be useful
in this direction:

Let $D_m$ be the $m$-dimensional Euclidean unit ball, and define
$d(t,m)=|tD_m\cap {\bf Z}^m|$ be the cardinality of the set of
integer points in $tD_m$. A simple lower bound for $d(t,m)$
can be given by counting the points with coordinates $0,\pm 1$ in $tD_m$:
$$d(t,m)\geq \sum_{k=0}^{[t^2]} \binom{n}{k} 2^k\geq \binom{n}{[t^2]}
2^{[t^2]}.\leqno (6.2)$$
By Theorem 2.4, for every $m\leq c_1\sqrt{n}$ we can find
$\sigma\subseteq\{1,\ldots,n\}$ with $|\sigma |=m$ and
$E\cap {\mathbf R}^{\sigma }\supseteq\frac{c_2}{\sqrt{m}M_E}D_{\sigma }$,
where $c_1, c_2>0$ are absolute constants. Assuming that $M_E<c_2$
and using (6.2) we have some non-trivial information: It is clear
that
$$|E\cap {\bf Z}^n|\geq \max_m\{|E\cap {\bf Z}^{\sigma }|:
|\sigma |=m\leq c_1\sqrt{n}\}\leqno (6.3)$$

\noindent Thus, we have:

\medskip

\noindent {\bf Fact II}: {\it Let $E$ be an ellipsoid in
${\mathbf R}^n$ with $M_E<c_2<1$. Then,}
$$|E\cap {\bf Z}^n|\geq
\max_m\{d(\frac{c_2}{\sqrt{m}M_E},m): m\leq c_1\sqrt{n}\}$$
$$\geq\max_m\{\binom{n}{[c_2^2/mM_E^2]}2^{[c_2^2/mM_E^2]}:m\leq c_1\sqrt{n}\}.$$

\medskip

The question of computation of the number of integer points inside
an ellipsoid (or, more generally, inside a symmetric convex body)
in ${\mathbf R}^n$ was relaxed in several directions in [M6]. One of
the questions stated asks for ``almost integer" points inside $E$
in the following precise sense: for a given $\theta\in (0,1)$,
find a projection of $E$ onto some coordinate subspace ${\mathbf R}^{\sigma }$
with $|\sigma |\geq (1-\theta )n$, which contains as many as possible
integer points. Then, $E$ itself will contain many points with
$[(1-\theta )n]$ coordinates which are distinct $[(1-\theta )n]$-dimensional
integers.

Our low $M^{\ast }$-estimate for ellipsoids provides an answer
to this question in terms of $M_E$. We know that there exists
$\sigma\subseteq\{1,\ldots,n\},\;|\sigma |= [(1-\theta )n]$,
such that
$${\mathcal P}_{\sigma }(E)\supseteq\frac{c\sqrt{\theta }}{\sqrt{\log (\frac
{2}{\theta })}M_E}\;D_{\sigma }.$$
This, and (6.2), lead to the following:

\medskip

\noindent {\bf Fact III}: {\it Let $E$ be an ellipsoid in
${\mathbf R}^n$. For every $\theta\in (0,1)$ there exists
$\sigma\subseteq\{1,\ldots,n\}$ with $|\sigma |=[(1-\theta )n]$
for which}
$$|{\mathcal P}_{\sigma }(E)\cap {\bf Z}^{\sigma }|\geq
d(\frac{c\sqrt{\theta }}{\sqrt{\log (\frac{2}{\theta })}M_E},[(1-\theta )n])$$
$$\geq \binom{[(1-\theta )n]}{[\frac{c^2\theta}{\log (\frac{2}{\theta })M_E^2}]}
2^{[\frac{c^2\theta }{\log (\frac{2}{\theta })M_E^2}]}.$$

Clearly, the results in Sections 3 and 4 give analogous estimates
for an arbitrary symmetric convex body.

\end{document}